\documentclass[11pt]{article}
\usepackage{amscd}
\usepackage{amsfonts}
\usepackage{amsmath}
\usepackage{amssymb}
\usepackage{amsthm}
\usepackage{bbm}
\usepackage{CJK}
\usepackage{fancyhdr}
\usepackage{graphicx}
\usepackage{hyperref}
\usepackage{indentfirst}
\usepackage{latexsym}
\usepackage{mathrsfs}
\usepackage{xypic}

\newtheorem{thm}{Theorem}[section]
\newtheorem{lemma}[thm]{Lemma}

\newcommand{\prf}{\mdseries \textbf{\textsc{Proof.}}\quad}

\usepackage[top=1in,bottom=1in,left=1.25in,right=1.25in]{geometry}
\def\<{\langle}
\def\>{\rangle}
\def\a{\alpha}
\def\b{\beta}

\def\c{\cdot}

\def\D{\Delta}

\def\g{\gamma}

\def\o{\otimes}

\def\v{\varepsilon}

\date{}
\begin{document}
\renewcommand{\baselinestretch}{1.2}
\renewcommand{\arraystretch}{1.0}
\title{\bf Enveloping actions  for  twisted partial actions  }
\author {{\bf Shuangjian Guo$^{1}$\footnote
        {Correspondence: Shuangjian Guo, School of Mathematics and Statistics, Guizhou University of Finance and Economics, Guiyang   550025, P. R. China;
        E-mail: shuangjguo@gmail.com}, Shengxiang Wang$^{2}$}\\
{\small 1. School of Mathematics and Statistics, Guizhou University of Finance and Economics}\\
{\small Guiyang  550025, P. R. China}\\
{\small 2. School of Mathematics and Finance, Chuzhou University}\\
 {\small Chuzhou 239000, P. R. China}}
 \maketitle
\begin{center}
\begin{minipage}{12.cm}

{\bf Abstract} Let $A\#_{\alpha, \omega}H$ be a partial crossed product.  In this paper, we first generalize the theorem about
the existence of an enveloping action to twisted  partial  actions. Second, we construct a Morita context between the partial crossed product and the crossed product related to the enveloping action. Furthermore, we discuss equivalences of partial crossed products  Finally,  we  investigate when $A\subset A\#_{\alpha, \omega}H$ becomes a
separable extension.
\smallskip

 {\bf Key words} Enveloping action, partial crossed product, Morita context, separable
extension.
\smallskip

 {\bf AMS(2010) Subject Classification}  16T05; 16T15
 \end{minipage}
 \end{center}
 \normalsize

\section{ Introduction}
\def\theequation{0. \arabic{equation}}
\setcounter{equation} {0}

Partial actions of groups have been introduced in the theory of operator algebras
as a general approach to study $C^{\ast}$-algebras by partial isometries in \cite{E94}.
A treatment from a purely algebraic point of view was given recently in \cite{D07}-\cite{DE10}. In particular, the Galois theory of partial group actions developed in
\cite{D07} inspired further Galois theoretic results in \cite{CD05}, as well as the introduction and
study of partial Hopf actions and coactions in \cite{C08}. The latter paper became in turn
the starting point for further investigation of partial Hopf (co)actions in \cite{A09}-\cite{A14}.
The general notion of a (continuous) twisted partial action of a locally compact
group on a $C^{\ast}$-algebra (a twisted partial $C^{\ast}$-dynamical system) and the corresponding crossed products were given by Exel in \cite{E97}. Twisted partial actions of groups on abstract rings and corresponding crossed products were recently introduced in
\cite{DE08}. More recent algebraic results on twisted partial actions and corresponding
crossed products were obtained in \cite{DE10} and \cite{B10}.  The algebraic concept of twisted partial actions also motivated the study
of projective partial group representations, the corresponding partial Schur Multiplier and
the relation to partial group actions with K-valued twistings in \cite{DN10} and \cite{DN12}, contributing
towards the elaboration of a background for a general cohomological theory based on partial
actions. Further information around partial actions may be consulted in the survey \cite{D11}.    As a unified approach for twisted
partial group actions, partial Hopf actions and twisted actions of Hopf algebras,
the notion of a twisted partial Hopf action is introduced in \cite{A14}.

Certainly, the theory of partial actions of Hopf algebras remains as a huge
landscape to be explored. Alves and  Batista \cite{A10}  generalized the theorem about
the existence of an enveloping action, also known as the globalization theorem and
they constructd a Morita context between the partial smash product and the smash
product related to the enveloping action, and this present work intends to generalize some results
for partial group actions and partial Hopf actions to the context of twisted partial  actions.

  This paper is organized as follows:

In Section 3,  we prove the  existence of an enveloping action for such a
twisted partial  actions. In Section 4,  we construct a Morita context between the partial crossed
 product $A\#_{\alpha,\omega}H$ and the crossed product $B\#_u H$,
where $H$ is a Hopf algebra which acts partially on the unital
algebra $A$, $B$ is an enveloping action for twisted partial actions. This
result can also be found in \cite{A10} for the context of partial group actions. In Section 5,  we shall discuss equivalences of partial crossed products and this result is similar to \cite{A14}. In Section 6,  we  investigate when $A\subset A\#_{\alpha, \omega}H$ becomes a
separable extension.

 \section{Preliminaries}
\def\theequation{\arabic{section}.\arabic{equation}}
\setcounter{equation} {0}

Throughout the paper, let $k$ be a fixed
 field and all algebraic systems are supposed to be over $k$.
 Let $M$ be a vector space
 over $k$ and let $id_M$  the usual identity map.
  For the comultiplication
 $\D $ in a coalgebra $C$ with a counit $\v _C$,
 we use the Sweedler-Heyneman's notation (see Sweedler \cite{S69}):
 $ \Delta(c)=c_{1}\o c_{2}$, for any $c\in C$.
 \vskip0.2cm

 We first recall some basic results and propositions that we will need later from Alves and  Batista \cite{A09},\cite{A10}.
 \medskip

  \noindent{\bf 2.1.  Partial  module algebra}~ Let $H$ be a Hopf algerba  and $A$ an  algebra. $A$ is said to be a partial
    $H$-module  algebra if there exists a $k$-linear map $\c\{\c: H\otimes A\rightarrow A\}$ satisfying the following conditions:
\begin{eqnarray*}
&&h\c (ab) =(h_{(1)})\c a)(h_{(2)}\c b),\\
&& 1_H\c a = a,\\
&&  h\c(g \c a) =(h_{(1)}\c 1_A)(h_{(2)}g\c a),
\end{eqnarray*}
for all $h,g\in H$ and $a, b\in A$.

\section{Enveloping actions}
\def\theequation{\arabic{section}.\arabic{equation}}
\setcounter{equation} {0}

Recall from \cite{A14} that  the definition of the partial crossed product. A twisted partial
action of $H$ on $A$ is a pair $(\alpha, \omega)$, where $\a: H\o A\rightarrow A, \a(h\o a)=h\c a$ and $\omega: H\o H\rightarrow A, \omega(h\o g)=\omega(h, g)$ be two linear maps such that the following
conditions hold:
\begin{eqnarray}
1_H\c a=a,\\
h\c(ab)=(h_{(1)}\c a)(h_{(2)}\c b),\\
(h_{(1)}\c (l_{(1)}\c a))\omega(h_{(2)}, l_{(2)})=\omega(h_{(1)}, l_{(1)})(h_{(2)}l_{(2)}\c a),\\
\omega(h, l)=\omega(h_{(1)}, l_{(1)})(h_{(2)}l_{(2)}\c 1_A),
\end{eqnarray}
for all $a,b\in A$ and $h, l\in A$.

{\bf Lemma 3.1.}$^{\cite{A14}}$ If $(\alpha, \omega)$ is a twisted partial action, then the following identities
hold:
\begin{eqnarray}
\omega(h, l)=(h_{(1)}\c (l_{(1)}\c 1_A))\omega(h_{(2)}, l_{(2)})=(h_{(1)}\c 1_A)\omega(h_{(2)}, l),
\end{eqnarray}
for all $h, l\in H$.

We say that the map $\omega$ is trivial, if the following condition holds:
\begin{eqnarray*}
h\c (l \c 1_A)=\omega(h, l)=(h_{(1)}\c 1_A)(h_{(2)}l\c 1_A),
\end{eqnarray*}
for all $h, l\in H$. In this case, the twisted partial action $(\alpha, \omega)$ turns out to be a
partial action of $H$ on $A$.

Given any twisted partial action $(\alpha, \omega)$ of $H$ on $A$, we can define on the vector
space $A\o H$ a product, given by the multiplication
\begin{eqnarray*}
(a\o h)(b\o l)=a(h_{(1)}\c b)\omega(h_{(2)}, l_{(1)})\o h_{(3)}l_{(2)},
\end{eqnarray*}
for all $a, b\in A$ and $h, l\in H$. Denote $A\#_{\a,\omega} H$ to be the subspace of $A\o H$ generated
by the elements of the form $a\# h=a(h_{(1)}\c 1_A)\o h_{(2)}$, for all $a\in A$ and $h\in H$.  Recall from [] that $A\#_{\a,\omega} H$ is associative via the multiplication Eq.(3.7)
and $A\#_{\a,\omega} H$ is unital with $1_A\# 1_H=1_A\o 1_H$  if and only if
\begin{eqnarray}
\omega(1_H, h)=\omega(h, 1_H)=h\c 1_A,\\
(h_{(1)}\c (l_{(1)}\c a))\omega(h_{(2)}, l_{(2)})=\omega(h_{(1)}, l_{(1)})(h_{(2)}l_{(2)}\c a),\\
(h_{(1)}\c \omega(l_{(1)}, m_{(1)}))\omega(h_{(2)}, l_{(2)}m_{(2)})=\omega(h_{(1)}, l_{(1)})\omega(h_{(2)}l_{(2)}, m_{(2)}),
\end{eqnarray}
for all $h,l,m\in H$ and $a\in A$.

Recall from \cite{A14} that let $B$ be a unital $k$-algebra measured by an action $\beta : H \o B\rightarrow B$, denoted by $\b(h\o b)=h\triangleright b$, which is twisted by a map $u: H\o H\rightarrow B$, i.e.,
\begin{eqnarray*}
1_H\triangleright a = a, \\
 h \triangleright(ab) = (h_{(1)} \triangleright
a)(h_{(2)}\triangleright b),\\
 h \triangleright 1_{A} = \varepsilon(h)1_{A},\\
u(1_{H}, h)  = u(h,
1_{H})= \varepsilon(h)1_{A},\\
(h_{(1)} \triangleright(l_{(1)} \triangleright a))u(h_{(2)},l_{(2)})
 =u(h_{(1)}, l_{(1)})(h_{(2)}l_{(2)}\triangleright a),\\
u(h_{(1)}, l_{(1)})u(h_{(2)}l_{(2)}, m) =
 (h_{(1)} \triangleright u(l_{(1)}, m_{(1)}))
u(h_{(2)}, l_{(2)}m_{(2)}),
\end{eqnarray*}
for all $h, l, m \in H$ and $a\in A$.  Suppose that $1_A$ is a non-trivial central idempotent of $B$, and let $A$ be the ideal generated
by $1_A$. Given $a\in A$ and $h\in H$, define a map $\c : H\o A\rightarrow A $ by $h\c a=1_A(h\triangleright a)$.
We say that the partial action $h\c a=1_A(h\triangleright a)$ and $\omega(h, l)=(h_{(1)}\c 1_A)u(h_{(2)}, l_{(1)})(h_{(3)}l_{(2)}\c 1_A)$ is the twisted partial action induced by $B$.

{\bf Definition 3.2.}$^{\cite{A14}}$
With notations as above. A morphism of algebras $\theta:A\rightarrow B$ is said to be a  morphism of  partial $H$-module algebras if
 $\theta(h\c a)=h\c \theta(a)$ for all $h,k\in H$ and
 $a\in A$.  If, in addition,  $\theta$ is an isomorphism, the partial actions are called equivalent.

Recall from \cite{A10} that if $B$ is an $H$-module algebra and $A$ is a right ideal of $B$ with
unity $1_A$, the induced partial action on $A$ is called admissible if $B=H\rhd A$.

{\bf Definition 3.3.}
With notations as above. An enveloping action for $A$ is a pair $(B,\theta),$  where

(a) $B$ is a unital $k$-algebra measured by an action $\beta : H \o B\rightarrow B$;

(b) The map $\theta:A\rightarrow B$ is a monomorphism of algebras;

(c) The sub-algebra $\theta(A)$ is an ideal in $B$;

(d) The partial action on $A$ is equivalent to the induced partial action on $\theta(A)$;

(e) The induced partial action on $\theta(A)$ is admissible.

(d) $\theta(a\omega(g, h))=\theta(a)u(g, h)$, $\theta(\omega(g, h)a)=u(g, h)\theta(a)$, for any $g, h\in H$ and $a\in A$.

\section{A Morita context}
\def\theequation{\arabic{section}.\arabic{equation}}
\setcounter{equation} {0}

In this section, we will construct a Morita context between the partial crossed product $A\#_{\alpha, \omega} H$ and the crossed product $B\#_{u} H$, where $B$ is an enveloping action for the partial crossed product.
\begin{lemma}
Let $(\alpha, \omega)$ be a twisted partial action and $(B,\theta)$ an enveloping action, then there is an algebra monomorphism from the
partial crossed product $A\#_{\alpha, \omega} H$ into the crossed product $B\#_{u} H$.
 \end{lemma}
\prf Define $\Phi:A\#_{\alpha, \omega} H\rightarrow B\#_{u} H$ by $a\otimes h \mapsto \theta(a)\o h$ for $h,g\in H$ and $a,b\in A$.
We first check that $\Phi$ is a morphism of algebras as follows:
\begin{eqnarray*}
\Phi((a\otimes h)(b\otimes g))
&=&\Phi( a(h_{(1)}\c b) \omega(h_{(2)}, g_{(1)})\otimes
h_{(3)}g_{(2)}).  \\
&=& \theta(a(h_{(1)}\c b) \omega(h_{(2)}, g_{(1)}))\#
h_{(3)}g_{(2)}\\
&=&  \theta(a)(h_{(1)}\c \theta(b))u(h_{(2)}, g_{(1)})\#
h_{(3)}g_{(2)}\\
&=& (\theta(a)\# h)(\theta(b)\# g)\\
&=& \Phi(a\otimes h)\Phi(b\otimes g).
\end{eqnarray*}

Next, we will verify that $\Phi$ is injective. For this purpose, take $x=\sum_{i=1}^{n}a_{i}\otimes h_{i} \in ker\Phi$ and choose $\{a_i\}^{n}_{i=1}$ to be linearly independent. Since $\theta$ is injective, we conclude that $\theta(a_{i})$ are linearly independent.
For each $f\in H^{\ast}$, $\sum_{i=1}^{n}\theta(a_{i})f( h_{i})=0,$  it follows that $f(h_{i})=0,$  so $ h_{i}=0$.
 Therefore we have $x=0$ and $\Phi$ is injective, as desired.

Since the partial crossed product $A\#_{\alpha, \omega} H$ is a subalgebra of $A\otimes H$, it is
 injectively mapped into $B\#_{u} H$ by $\Phi.$  A typical element of the image of the partial crossed product is
\begin{eqnarray*}
\Phi((a\otimes h)(1_A\otimes 1_{H}))
&=& \Phi(a\otimes h)\Phi(1_A\otimes 1_{H})\\
&=& (\theta(a)\# h)(\theta(1_A)\# 1_{H})\\
&=& \theta(a(h_{(1)}\c \theta(1_A)))\#
h_{(2)}.
\end{eqnarray*}
And this completes the proof. $\hfill \Box$
\medskip

Take $M=\Phi(A\otimes H)=\{\sum_{i=1}^{n}\theta(a_{i})\#  h_{i}; a_{i}\in A\}$  and take $N$ as the subspace of
$B\#_{u} H$  generated by the elements $h_{(1)}\c \theta(a)\otimes
h_{(2)}$ with $h\in H$ and $a\in A.$

{\bf Proposition 4.2.}
Let $(\alpha, \omega)$ be a twisted partial action and $(B,\theta)$ an enveloping action.  Suppose that $\theta(A)$ is an ideal of $B$,  then $M$ is a right $B\#_u H$
module and $N$ is a left $B\#_u H$ module.

{\bf Proof.} In order to prove $M$ is a right $B\#_u H$ module, let $\theta(a)\# h\in M$ and $b\# k\in B\#_u H$,  then
              $$(\theta(a)\# h)(b\# k)= \theta(a)(h_{(1)}\triangleright b)u(h_{(2)}, k_{(1)} )\# h_{(3)}k_{(2)}. $$
Which lies in $\Phi(A\otimes H)$ because $\theta(A)$ is an ideal in $B$.

Now we show that $N$ is a left $B\# H$ module is similar to \cite{A10}.
$\hfill \Box$
\medskip

Now, the left $A\#_{\alpha, \omega} H$ module structure on $M$ and a right $A\#_{\alpha, \omega} H$ module
structure on $N$ induced by the monomorphism $\Phi$ is in \cite{A10}.

{\bf Proposition 4.3.}
Under the same hypotheses of Proposition 4.2, $M$ is indeed a
left $A\#_{\alpha, \omega} H$  module with the map $\blacktriangleright$ and N is a right $A\#_{\alpha, \omega} H$  module with the map $\blacktriangleleft$.

{\bf Proof.} We first claim that $A\#_{\alpha, \omega} H\blacktriangleright M \subseteq M.$ In fact,
\begin{eqnarray*}
&&(a( h_{(1)}\c 1_A)\otimes h_{(2)})\blacktriangleright (\theta(b))\# k)\\
&=&
(\theta(a)(h_{(1)}\c \theta(1_A))\# h_{(2)})(\theta(b))\# k)\\
&=& \theta(a)(h_{(1)}\c \theta(1_A))(h_{(2)}\triangleright \theta(b))u(h_{(3)}, k_{(1)})\# h_{(4)}k_{(2)}\\
&=& \theta(a)(h_{(1)}\c \theta(1_A))(h_{(2)}\c \theta(b))u(h_{(3)}, k_{(1)})\# h_{(4)}k_{(2)}\\
&=& \theta(a)(h_{(1)}\c \theta(b)) u(h_{(2)}, k_{(1)})\# h_{(3)}k_{(2)}.
\end{eqnarray*}
Which lies inside $M$ because $\theta(A)$ is an ideal of $B$.
$\hfill \Box$

Next, we verify that $N \blacktriangleleft A\#_{\alpha, \omega} H\subseteq N$,  which is similarly to  $N$ is a left $B\#_u H$ module.
which holds because  $\theta(1_A)$ is  a central idempotent.

The last ingredient for a Morita context is to define two bimodule
morphisms
$$\sigma: N\otimes_{A\#_{\alpha, \omega} H} M \rightarrow B\#_u H ~~\mbox{and }$$
$$\tau: M\otimes_{B\#_u H} N\rightarrow A\#_{\alpha, \omega} H\cong \Phi(A\#_{\alpha, \omega} H).$$
As $M, N$ and $A\#_{\alpha, \omega} H$ are viewed as subalgebras of $B\# H,$  these two maps can be taken
as the usual multiplication on $B\#_u H.$
The associativity of the product assures us
that these maps are bimodule morphisms and satisfy the associativity conditions.

{\bf Proposition 4.4.}
 The partial crossed
 product $A\#_{\alpha, \omega} H$ is Morita equivalent to the  crossed product $B\#_u H$.

\section{Equivalences of partial crossed products}
\def\theequation{\arabic{section}.\arabic{equation}}
\setcounter{equation} {0}

In this section, we shall discuss equivalences of partial crossed products and this result is similar to \cite{A14}. From now on, unless explicitly stated, we always assume partial actions of a Hopf algebra $H$ over an algebra $A$ such that the map $e\in$ Hom$(H, A)$ ,
given by $e(h)=h\c 1_A$, is central with respect to the convolution product. Given a twisted partial H-module algebra $A=(A, \alpha, \omega)$, we define two linear maps $f_1, f_2 : H\o H \rightarrow A$ as follows:
\begin{eqnarray*}
f_1(h, k)=(h\c 1_A)\varepsilon(k),~~~~~f_2(h, k)=hk\c 1_A.
\end{eqnarray*}

{\bf Definition 5.1.}$^{\cite{A14}}$ Let $A=(A, \alpha, \omega)$ be a twisted partial $H$-module algebra. We will
say that the partial action is symmetric, if\\
(1) $f_1, f_2$  are central in Hom ($H\o H, A$),\\
(2) $\omega$ satisfies the conditions Eq.(3.7) and Eq.(3.8) and has a convolution inverse $\omega'$
in the ideal $<f_1\ast f_2>\in$ Hom ($H\o H, A$),\\
(3) $ h\c(k\c 1_A)=(h_{(1)} \c 1_A)(h_{(2)}k\c 1_A)$ for any $h, k \in H$.

{\bf Definition 5.2.} If there exists linear maps $u, v\in$ Hom$(H,A)$ satisfy $(u\ast v)(h)=(v\ast u)(h)=h\c 1_A$, $u(h)=u(h_{(1)})(h_{(2)}\c 1_A)=(h_{(1)}\c 1_A)u(h_{(2)})$ and $u(1_H)=v(1_H)=1_A$. Then we call $v$ is a weak convolution-invertible linear map and $u=v^{-1}$.

 Let $H$ be a Hopf algebra and  $(A, \alpha, \omega)$  a symmetric twisted partial $H$-module algebra, and
$v\in$ Hom$(H,A)$ a weak convolution-invertible linear map.  Define
$\omega^{v}:H\otimes H\rightarrow A $ and weakly action of $H$ on
$A$  by
$$
\omega^{v}(h,g) =
v(h_{(1)})(h_{(2)}\c
v(g_{(1)}))\omega(h_{(3)},g_{(2)})
v^{-1}(h_{(4)}g_{(3)})
$$
and
$$
 h\c^v a =
v(h_{(1)})(h_{(2)}\c
a)v^{-1}(h_{(3)})
$$
for any $h, g\in H$ and $a\in A$.

{\bf Lemma 5.3.} With the notations as above, let $H$ be a Hopf algebra and  $(A, \alpha, \omega)$  a symmetric twisted partial $H$-module algebra. Then
 $\omega^{vu}=(\omega^{u})^{v}$
 and $\c^{vu}=(\c^{u})^{v}$, where $u, v\in$ Hom$(H,A)$ are weak convolution-invertible linear maps.

{\bf Proof.} For any $h, g\in H, a\in A$, we have
we have
\begin{eqnarray*}
&&\omega^{vu}(h,g)\\
&=&  (vu(h_{(1)})(h_{(2)}\c
(vu(g_{(1)})))\omega(h_{(3)},g_{(2)})
(vu)^{-1}(h_{(4)}g_{(3)})\\
&=&  v(h_{(1)})u( h_{(2)})(h_{(3)}\c v(g_{(1)}))(h_{(4)}\c u(g_{(2)})\omega(h_{(5)},g_{(3)})u^{-1}(h_{(6)}g_{(4)})v^{-1}(h_{(7)}g_{(5)})\\
&=& v(h_{(1)})u( h_{(2)})(h_{(3,1)}\c v(g_{(1)}))u^{-1}(h_{(4)})u( h_{(5)})(h_{(6)}\c u(g_{(2)}))\\
&&\omega(h_{(7)},g_{(3)})u^{-1}(h_{(8)}g_{(4)})v^{-1}(h_{(9)}g_{(5)})\\
&=&v(h_{(1)})(h_{(2)}\c^u v(g_{(1)}))\omega^{u}(h_{(3)},g_{(2)})v^{-1}(h_{(4)}g_{(3)})\\
&=&(\chi^{u})^{v}(h,g)
\end{eqnarray*}
and thus
 $\chi^{vu}=(\chi^{u})^{v}$.

Also,
\begin{eqnarray*}
  && h(\c^{u})^{v}a\\
&=& v(h_{(1)})(h_{(2)}\c^{u}a)v^{-1}(h_{(3)})\\
&=&  v(h_{(1)})u(h_{(2)})(h_{(3)}\c a)u^{-1}(h_{(4)})v^{-1}(h_{(5)})\\
&=&  v(h_{(1)})u(h_{(2)}))(h_{(3)}\c a) u^{-1}(h_{(4)})v^{-1}(h_{(5)})\\
&=& h\c^{vu}a
\end{eqnarray*}
and so
$\c^{vu}=(\c^u)^{v}$

 This completes the proof. $\hfill \Box$

{\bf Theorem 5.4.}  Let $H$ be a Hopf algebra and $(A, \alpha, \omega)$  a symmetric twisted partial $H$-module algebra, and
 $v\in$ Hom$(H,A)$ a
 convolution-invertible linear map, with the above notations $\omega^{v}, \c^v$.
 Then we have the following assertions:

(1) As algebras, $A\#_{\a, \omega} H \cong
 A\#_{\a, \omega^{v}}H $;

(2) $\omega$ satisfies Eq.(3.6) if and only if  $\omega^{v}$
 satisfies Eq.(3.6);

(3) $(\omega,\c)$ satisfies Eq.(3.7) if and only if
 $(\omega^{v}, \c^{v})$ satisfies Eq.(3.7);

(4) If $(\omega,\c)$ satisfies Eq.(3.7), then
 $(\omega,\c)$
  satisfies Eq.(3.8)  if and only if $(\omega^{v}, \c^{v})$ satisfies Eq.(3.8);

(5) $A\#_{\a, \omega} H $ is a partial crossed product algebra if
 and only if $A\#_{\a, \omega^{v}}H $ is a
partial crossed product algebra , and they are isomorphic.

\prf (1) Define
$\Phi:A\#_{\a, \omega}H\mapsto  A\#_{\a, \omega^{v}}H$ by
 $ a\# h\rightarrow  av(h_{(1)})\# h_{(2)}$, for any $h,g\in H$, we have
\begin{eqnarray*}
&&\Phi((a\# h)(b\# g))\\
&=&\Phi( a(h_{(1)}\c ^{v}b)\omega^{v}(h_{(2)},g_{(1)})\# h_{(3)}g_{(2)})\\
&=& a(h_{(1)}\c^{v}b)\omega^{v}(h_{(2)},g_{(1)})v(h_{(3)}g_{(2)})\#  h_{(4)}g_{(3)})\\
&=&  av(h_{(1)})(h_{(2)}\c b)v^{-1}h_{(3)}v(h_{(4)})h_{(5)}\c \omega(h_{(6)},g_{(2)})v^{-1}(h_{(7)}g_{(3)})
v(h_{(8)}g_{(4)})\# h_{(9)}g_{(5)}\\
&=&  av(h_{(1)})(h_{(2)}\c b)(h_{(3)}\c v(g_{(1)})\omega (h_{(4)},g_{(2)})
\# h_{(5)}g_{(3)}\\
&=& \Phi(a\# h)\Phi(b\# g).
\end{eqnarray*}
Clearly $\Phi$ is bijective, $\Phi^{-1}(a\# h)=\sum av^{-1}(h_{(1)})\# h_{(2)}$
$a,b\in A,h\in H$, since
\begin{eqnarray*}
&&\Phi\Psi(a\# h)\\
&=&\Phi( av^{-1}(h_{(1)})\# h_{(2)})\\
&=& av^{-1}(h_{(1)})v(h_{(2)})\# h_{(3)}\\
&=& av^{-1}(h_{(1)})v(h_{(2)})\# h_{(3)} \\
&=& a\# h.
\end{eqnarray*}

(2) Straightforward.

 (3) If $(\omega,\c)$ satisfies
 Eq.(3.7), then
\begin{eqnarray*}
&& (h_{(1)}\c^{v}(g_{(1)}\c^{v}(a)))\omega^{v}(h_{(2)},g_{(2)})\\
&=&  v(h_{(1)})(h_{(2)}\c (v(g_{(1)}))(g_{(2)}\c a)v^{-1}(g_{(3)}) v^{-1}(h_{(3)})
v(h_{(4)})(h_{(5)}\c v(g_{(4)})\\
&&\hspace{7cm}\omega(h_{(6)},g_{(5)})v^{-1}(h_{(7)}g_{(6)})\\
&=&  v(h_{(1)})(h_{(2)}\c v(g_{(1)}))(h_{(3)}\c (g_{(2)}\c a))\omega(h_{(4)},g_{(3)})v^{-1}(h_{(5)}g_{(4)})\\
&=&  v(h_{(1)})(h_{(2)}\c v(g_{(1)}))\omega(h_{(3)},g_{(2)})(h_{(4)}g_{(3)}\c a)v^{-1}(h_{(5)}g_{(4)})\\
&=& \omega^{v}(h_{(1)},g_{(1)})(h_{(2)}g_{(2)}\c^{v} a).
\end{eqnarray*}
Conversely, we get it from Lemma 5.3.

 (4) If $(\omega,\c)$  satisfies Eq.(3.7)and Eq.(3.8), then for $h, g, m\in
 H$, we have
\begin{eqnarray*}
&& (h_{(1)}\c^{v}\omega^{v}(g_{(1)},m_{(1)}))\omega^{v}(h_{(2)},g_{(2)}m_{(2)})\\
&=&  v(h_{(1)})(h_{(2)}\c[(v (g_{(1)})(g_{(2)}\c^{v}(m_{(1)})))
\omega(g_{(3)},m_{(2)})v^{-1}(g_{(4)}m_{(3)})])\\
&&
v^{-1}(h_{(3)})v(h_{(4)})
(h_{(5)}\c v(g_{(5)}m_{(4)}))\omega(h_{(6)},g_{(6)}m_{(5)})
v^{-1}(h_{(7)}g_{(7)}m_{(6)}\\
&=&  v(h_{(1)})(h_{(2)}\c v (g_{(1)}))(h_{(3)}\c (g_{(2)}\c v (m_{(1)})))\\
&&(h_{(4)}\c \omega(g_{(3)}m_{(2)}))\omega(h_{(5)},g_{(4)}m_{(3,1)})v^{-1}(h_{(6)}g_{(5)}m_{(4)})\\
&=&  v(h_{(1)})(h_{(2)}\c v (g_{(1)}))(h_{(3)}\c (g_{(2)}\c v (m_{(1)})))\\
&&\omega (h_{(4)},g_{(3)})\chi(h_{(5,1)},g_{(4,1)}m_{(2,1)})v^{-1}(h_{(6)}g_{(5)}m_{(4)})\\
&=&  v(h_{(1,\alpha)})(h_{(2)}\c v (g_{(1)}))\omega(h_{(3)}g_{(2)})(h_{(4)}g_{(3)}\c v (m_{(1)}))\\
&&\omega(h_{(5)}g_{(4)},m_{(2)})v^{-1}(h_{(6)}g_{(5)}m_{(4)})\\
&=&  \omega^{v}(h_{(1)},g_{(1)})\omega(h_{(2)}g_{(2)},m).
\end{eqnarray*}
Conversely, we get it from Lemma 5.3.

 (5) Clearly. $\hfill \Box$

\section{Separable extension for partial crossed products}
\def\theequation{\arabic{section}.\arabic{equation}}
\setcounter{equation} {0}

In this section, we will investigate when $A\subset A\#_{\alpha, \omega}H$ is a separable extension.

{\bf Definition 6.1.}$^{\cite{A14}}$ Let $B$ be a right $H$-comodule unital algebra with coaction given by $\rho : B \rightarrow
B\o H $ and let $A$ be a subalgebra of $B$. We will say that $A \subset B$ is an $H$-extension if $A = B^{coH}$.
An $H$-extension $A \subset B$ is partially cleft if there is a pair of $k$-linear maps
$\g, \g' : H \rightarrow B$
such that\\
(i) $\g(1_H)=1_B,$\\
(ii)$ \rho\circ \g=(\g\o id_H)\Delta$ and $\rho\circ \g'=(\g'\o S)\Delta^{cop}$,\\
(iii)$(\g\ast\g') \circ M$ is a central element in the convolution algebra Hom$(H \o H, A)$, where
$M : H \o H \rightarrow H $ is the multiplication in $H$, and $(\g\ast\g') (h)$ commutes with every
element of $A$ for each $h \in H$.

{\bf Lemma 6.2} Let $c\in C_{A\#_{\a, \omega}H}(A)$, the centralizer of $A$ in $A\#_{\a, \omega}H$, and assume that $H$ is
cocommutative. Then for $h \in H$,
\begin{eqnarray*}
\g'(h_{(1)})c(S(h_{(2)})\c 1_A)\g(h_{(3)})&=&\g(S(h_{(2)}))c\g'(S(h_{(1)}))\\
&=& S(h)\c c.
\end{eqnarray*}

{\bf Proof.} For any $h\in H$, we have
\begin{eqnarray*}
&&\g'(h_{(1)})c(S(h_{(2)})\c 1_A)\g(h_{(3)})\\
&=&\g'(h_{(1)})c\g(h_{(4)})\g(S(h_{(3)}))\g'(S(h_{(2)}))\\
&=&\g'(h_{(1)})c\g(h_{(2)})\g(S(h_{(3)}))\g'(S(h_{(4)}))\\
&=& \g'(h_{(1)})\g(h_{(2)})\g(S(h_{(3)}))c\g'(S(h_{(4)}))\\
&=& (S(h_{(1)})\c 1_A)\g(S(h_{(2)}))c\g'(S(h_{(3)}))\\
&=&(S(h_{(1)})\c 1_A)(S(h_{(2)})\c c)\\
&=&S(h)\c c
\end{eqnarray*}

{\bf Theorem 6.2.} Let $H$ be a finite dimensional cocommutative  Hopf algebra
and $(A, \c, (\omega, \omega'))$ is a symmetric partial twisted $H$-module algebra, and
$t\neq 0$ be a left integral in $H$. Assume that there exists $c$ in the center of $A$ such that
$t\c c=1_A$. Then $A\subset A\#_{\alpha, \omega}H$ is a separable extension.

{\bf Proof. } In order to show that $A\subset A\#_{\alpha, \omega}H$ is separable, it suffices to find a separability
idempotent $e \in (A\#_{\alpha, \omega}H) \o_A (A\#_{\alpha, \omega}H)$; this means that\\
(1) $((a\#_{\alpha, \omega}h) \o_A (1\#_{\alpha, \omega} 1))e = e((1\#_{\alpha, \omega} 1) \o_A (a\#_{\alpha, \omega}h))$ for all $a\#_{\alpha, \omega}h \in A\#_{\alpha, \omega}H$ and\\
(2) $m_{A\#_{\alpha, \omega}H}(e)= 1\#_{\alpha, \omega}1$, where $m_{A\#_{\alpha, \omega}H}$ denotes multiplication $(A\#_{\alpha, \omega}H) \o (A\#_{\alpha, \omega}H) \rightarrow
(A\#_{\alpha, \omega}H)$.

Now since $t$ is a left integral in $H$, $u=S(t)$ is a right integral in $H$. We claim that $f=\g'(u_{(1)})\o _A\g(u_{(2)})\in (A\#_{\alpha, \omega}H) \o_A (A\#_{\alpha, \omega}H)$  satisfies condition
(1) above. By \cite{C08}, one can get the canonical map can:
\begin{eqnarray*}
A\#_{\alpha, \omega}H) \o_A (A\#_{\alpha, \omega}H) \rightarrow (A\#_{\alpha, \omega}H)\o H\\
(a\#_{\alpha, \omega}h)(b\#_{\alpha, \omega}g)\mapsto (a\#_{\alpha, \omega}h)(b\#_{\alpha, \omega}g_{(1)})\o g_{(2)}
\end{eqnarray*}
is a bijection,  since by \cite{A14} $A\subset A\#_{\alpha, \omega}H$ is a partially cleft extension, we can check $A\subset A\#_{\alpha, \omega}H$ also is partial $H$-Hopf Galois, we claim that both $((a\#_{\alpha, \omega}h) \o_A (1\#_{\alpha, \omega} 1))f$ and $ f((1\#_{\alpha, \omega} 1) \o_A (a\#_{\alpha, \omega}h))$ go to $(a\#_{\alpha, \omega}h)(S(u_{(1)})\c 1_A) \o u_{(2)}$. As a matter of fact,
\begin{eqnarray*}
&&can(((a\#_{\alpha, \omega}h) \o_A (1\#_{\alpha, \omega} 1))f)\\
&=& can((a\#_{\alpha, \omega}h)\g'(u_{(1)})\o_A\g(u_{(2)}))\\
&=&(a\#_{\alpha, \omega}h)\g'(u_{(1)})\g(u_{(2)})\o u_{(3)}\\
&=&(a\#_{\alpha, \omega}h)(S(u_{(1)})\c 1_A)\o u_{(2)},
\end{eqnarray*}
and
\begin{eqnarray*}
&&can((f(a\#_{\alpha, \omega}h) \o_A (1\#_{\alpha, \omega} 1)))\\
&=&can(\g'(u_{(1)})\o_A\g(u_{(2)})(a\#_{\alpha, \omega}h))\\
&=&can(\g'(u_{(1)})\o_A (u_{(2)}\c a)\omega(u_{(3)}, h_{(1)})\g(u_{(4)}h_{(2)}))\\
&=& \g'(u_{(1)})(u_{(2)}\c a)\omega(u_{(3)}, h_{(1)})\g(u_{(4)}h_{(2)})\o u_{(5)}h_{(3)}\\
&=&\g'(u_{(1)})\g(u_{(2)}) a\g'(u_{(3)})\omega(u_{(4)}, h_{(1)})\g(u_{(5)}h_{(2)})\o u_{(6)}h_{(3)}\\
&=&(S(u_{(1)})\c 1_A)a\g'(u_{(2)})\omega(u_{(3)}, h_{(1)})\g(u_{(4)}h_{(2)})\o u_{(5)}h_{(3)}\\
&=&(S(u_{(1)})\c 1_A)a\g(h_{(1)})\g'(u_{(2)}h_{(2)})\g(u_{(3)}h_{(3)})\o u_{(4)}h_{(4)}\\
&=&(S(u_{(1)})\c 1_A)a\g(h_{(1)})\g'(u_{(2)})\g(u_{(3)})\o u_{(4)}\\
&=&(a\#_{\alpha, \omega}h)(a\#_{\alpha, \omega}h)(S(u_{(1)})\c 1_A)\o u_{(2)}.
\end{eqnarray*}

Now since $c$ is in the center of $A$, the map
\begin{eqnarray*}
(A\#_{\alpha, \omega}H) \o_A (A\#_{\alpha, \omega}H) \rightarrow (A\#_{\alpha, \omega}H) \o_A (A\#_{\alpha, \omega}H) \\
(a\#_{\alpha, \omega}h) \o_A (b\#_{\alpha, \omega}g) \mapsto (a\#_{\alpha, \omega}h)c \o_A (b\#_{\alpha, \omega}g)
\end{eqnarray*}
is well-defined. Thus $e=\g'(u_{(1)})c(S(u_{(2)})\c 1_A)\o _A\g(u_{(3)})$ also satisfies (1).

We also claim that $m(e) =\g'(u_{(1)})c(S(u_{(2)})\c 1_A)\g(u_{(3)})=1_A$. For, by Lemma 6.1,
\begin{eqnarray*}
\g'(u_{(1)})c(S(u_{(2)})\c 1_A)\g(u_{(3)})=S(u)\c c=S^2(t)\c c=t\c c=1_A.
\end{eqnarray*}
It follows that $A\subset A\#_{\alpha, \omega}H$ is a separable
extension.
\section*{\bf Acknowledgement}

 The work was  supported  by the TianYuan Special Funds of the National Natural Science Foundation of China (11426073), the NSF of Jiangsu Province (BK2012736) and the Fund of Science and Technology Department of Guizhou Province (2014GZ81365).

\end{document}